\documentclass[12pt]{article}
\usepackage{amssymb}
\newcommand{\dd}{{\rm \kern 3pt I\kern-9pt d}}

\newcommand{\Abar}{{\backslash\kern-8pt A}}

\title{\large ON SOME ERRORS RELATED TO THE GRADUATION OF MEASURING INSTRUMENTS}
\author{\sc Nicolas Bouleau\footnote{Ecole des Ponts, ParisTech, 28 rue des Saints P\`eres, 75007 Paris, France; 
e-mail : {\tt bouleau@enpc.fr}}}
\date{October 2006}

\font\helv Helvetica at 14pt
\newcommand{\CC}{\mbox{{\helv \i}{\rm\kern-5pt C}}}

\begin{document}
\maketitle

\noindent{\bf Abstract.} The error on a real quantity $Y$ due to the graduation of the measuring instrument may be approximately represented, when the graduation is regular and fines down, by a Dirichlet form on $\mathbb{R}$ whose square field operator does not depend on the probability law of $Y$ as soon as this law possesses a continuous density. This feature is related to the ``arbitrary functions principle" (Poincar\'e, Hopf). We give extensions of this property to $\mathbb{R}^d$ and to the Wiener space for some approximations of the Brownian motion and apply these results to the discretization of stochastic differential equations. \\


\noindent{\Large\textsf{I. Introduction.}}\\

\noindent{\large \textsf{I.1.}} The approximation of a random variable $Y$ by an other one $Y_n$ yields most often a Dirichlet form.  The framework is general, cf. Bouleau (2006) : Let $Y$ and $Y_n$ be defined on $(\Omega, \mathcal{A}, \mathbb{P})$ with values in the measurable space $(E,\mathcal{F})$, denoting $\mathbb{P}_Y$ the law of $Y$, if there exists an algebra $\mathcal{D}$ of bounded functions dense in $L^2(\mathbb{P}_Y)$ and a sequence of positive real numbers $\alpha_n$ such that, for all $\varphi\in \mathcal{D}$ there exists $\tilde{A}[\varphi]\in L^2(\mathbb{P}_Y)$ such that
$$\forall\psi\in\mathcal{D}\quad\alpha_n\mathbb{E}[(\varphi(Y_n)-\varphi(Y))(\psi(Y_n)-\psi(Y))]\rightarrow -2<\tilde{A}[\varphi],\psi>_{L^2(\mathbb{P}_Y)}$$
then $\mathcal{E}[\varphi,\psi]=-<\tilde{A}[\varphi],\psi>$ is a Dirichlet form.

{\it Often}, when this Dirichlet form exists and does not vanish, the conditional law of $Y_n$ given $Y=y$ is not reduced to a Dirac mass, and the variance of this conditional law yields the square field operator $\Gamma$. On the other hand when the approximation is deterministic, i.e. when $Y_n$ is a function of $Y$ say $Y_n=\eta_n(Y)$, then {\it most often} the symmetric bias operator $\tilde{A}$ and the Dirichlet form vanish, Bouleau (2006, examples 2.1 to 2.9 and remark 5).

Nevertheless, there are cases where the conditional law of $Y_n$ given $Y$ is a Dirac mass, i.e. $Y_n$ is a deterministic function of $Y$, and where the approximation of $Y$ by $Y_n$ yields even so a non zero Dirichlet form on $L^2(\mathbb{P}_Y)$.

This phenomenon is interesting, insofar as randomness (here the Dirichlet form) is generated by a deterministic device. In its simplest form, the phenomenon appears precisely when a quantity is measured by a graduated instrument when looking for the asymptotic limits as the graduation fines down. 

We will first expose (section {\textsf{I.2}}) the simplest case of a real quantity measured with an equidistant graduation. The mathematical argument is here the same as for the arbitrary functions method, about which we give a short historical survey (section {\textsf{I.3}}).

The second part is devoted to theoretical tools that we shall use later on. We first recall
the bias operators and the properties of the Dirichlet form associated with an approximation. Next we prove a version of a Girsanov theorem for Dirichlet forms which has its own interest, i.e. an answer to the question of an absolutely continuous change of measure for Dirichlet forms. At last we recall some simple properties of Rajchman measures.

Thanks to these tools, in the third part we take up again the classical case in order to extend it to $\mathbb{R}^d$ and to more general graduation.

The fourth part is concerned by Rajchman martingales, i.e. continuous local martingales whose brackets are a.s. Rajchman measures. These martingales possess remarkable properties of weak convergence for some approximations. Here the followed method is essentially an extension of the seminal idea of Rootz\'en (1980).
If we restrict the framework to the Wiener space, the obtained limit quadratic forms may be shown to be closable hence Dirichlet, this is done in the fifth part. This yields asymptotic speed of weak convergence for the discretization of stochastic differential equations in the case encountered for mechanical systems which was not explicited in the recent works of Kurtz and Protter (1991), Jacod and Protter (1998).\\

\noindent{\large \textsf{I.2.}} {\bf The basic example.}

Let $Y$ be a real random variable. It is approximated by $Y_n$ to the nearest graduation, i.e. 
$$Y_n=\frac{[nY]}{n}+\frac{1}{2n}.$$
($[x]$ denotes the entire part of $x$, and $\{x\}=x-[x]$ the fractional part).

We put $Y_n=Y+\xi_n(Y)$ where the function $\xi_n(x)=\frac{[nx]}{n}-\frac{1}{2n}-x$ is periodic with period $\frac{1}{n}$ and may be written $\xi_n(x) =\frac{1}{n}\theta(nx)$ with $\theta(x)=\frac{1}{2}-\{x\}.$

\setlength{\unitlength}{0.6pt}
\begin{picture}(700,300)(35,10)
\put(20,30){\line(1,0){270}}
\put(50,20){\line(0,1){250}}
\put(224,20){\line(0,1){250}}
\put(200,20){\line(0,1){250}}
\put(176,20){\line(0,1){250}}
\put(152,20){\line(0,1){250}}
\put(128,20){\line(0,1){250}}
\put(50,30){\line(1,1){220}}
\put(116,108){\circle*{3}}
\put(116,108){\line(1,0){24}}
\put(140,132){\circle*{3}}
\put(140,132){\line(1,0){24}}
\put(164,156){\circle*{3}}
\put(164,156){\line(1,0){24}}
\put(188,180){\circle*{3}}
\put(188,180){\line(1,0){24}}
\put(212,204){\circle*{3}}
\put(212,204){\line(1,0){24}}
\put(60,250){$Y_n$}
\put(30,15){O}\put(260,40){$Y$}
\put(370,130){\line(1,0){270}}
\put(400,60){\line(0,1){210}}
\put(574,60){\line(0,1){210}}
\put(550,60){\line(0,1){210}}
\put(526,60){\line(0,1){210}}
\put(502,60){\line(0,1){210}}
\put(478,60){\line(0,1){210}}
\put(466,142){\circle*{3}}
\put(466,142){\line(1,-1){24}}
\put(490,142){\circle*{3}}
\put(490,142){\line(1,-1){24}}
\put(514,142){\circle*{3}}
\put(514,142){\line(1,-1){24}}
\put(538,142){\circle*{3}}
\put(538,142){\line(1,-1){24}}
\put(562,142){\circle*{3}}
\put(562,142){\line(1,-1){24}}
\put(410,250){$\xi_n(x)$}
\put(380,112){O}
\put(630,140){$x$}
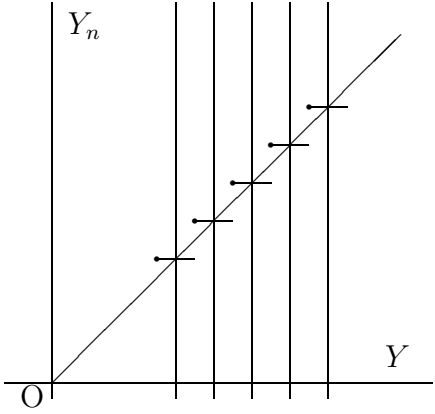
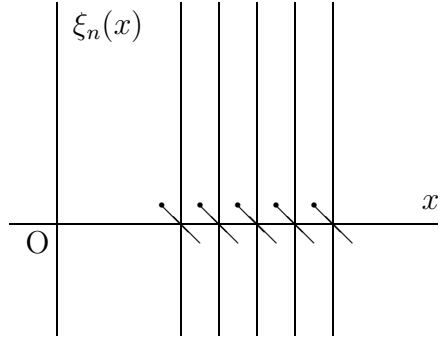
\put(330,30){\footnotesize fig 1 and 2: Approximation of the quantity Y} \put(390,10){\footnotesize to the nearest graduation.}
\end{picture}\\

Let $\mathbb{P}_Y$ the law of $Y$, we approximate Y by $Y_n$ on the algebra $\mathcal{D}=\mathcal{C}^1\cap Lip$ with the sequence $\alpha_n=n^2$, cf. Bouleau (2006).

Let us recall that the Rajchman class is the set of bounded measures on $\mathbb{R}$ whose Fourier transform vanishes at infinity. These measures are continuous (do not charge points) and are a band in the space of bounded measures on $\mathbb{R}$ (see section {\textsf{II.3}} below).\\

\noindent{\bf Theorem 1.} {\it If $\mathbb{P}_Y$ is a Rajchman measure,
\begin{equation}
(n(Y_n-Y),Y)\quad\stackrel{d}{\Longrightarrow}\quad(V,Y)
\end{equation}
where $V$ is uniform on $(-\frac{1}{2},\frac{1}{2})$ independent of $Y$, and for $\varphi\in\mathcal{C}^1\cap Lip$
\begin{equation}n^2\mathbb{E}[(\varphi(Y_n)-\varphi(Y))^2]\quad\rightarrow\quad\frac{1}{12}\mathbb{E}_Y[\varphi^{\prime 2}].\end{equation}}

\noindent Here $\stackrel{d}{\Longrightarrow}$ denotes the weak convergence, i.e. the convergence of probability measures on continuous bounded functions, $\mathbb{E}_Y$ is the expectation under $\mathbb{P}_Y$.

\noindent{\bf Proof.}  It is equivalent to study the weak convergence of  $(\frac{1}{2}+n(Y_n-Y),Y)=(\frac{1}{2}+\theta(nY),Y)$. Since $\frac{1}{2}+\theta$ takes its values in the unit interval, it is enough to study the convergence on the characters of the group $\mathbb{T}^1\times\mathbb{R}$, i.e. 
$$\mathbb{E}[e^{2i\pi k(\frac{1}{2}+\theta(nY))}e^{iuY}]=\mathbb{E}[e^{-2i\pi knY}e^{iuY}]=\Psi_Y(u-2\pi kn)$$
where $\Psi_Y$ is the characteristic function of $Y$. This tends to $\Psi(u)1_{\{k\neq0\}}$ since  $ \mathbb{P}_Y$ is Rajchman, proving the first assertion.

If $\varphi\in\mathcal{C}^1\cap Lip$, the relation $\varphi(Y_n)-\varphi(Y)=(Y_n-Y)\int_0^1\varphi^\prime(Y+\alpha(Y_n-Y))d\alpha$ gives
$$n^2\mathbb{E}[(\varphi(Y_n)-\varphi(Y))^2]=\mathbb{E}[\theta^2(nY)\varphi^{\prime 2}(Y)]+o(1)$$ and $\mathbb{E}[\theta^2(nY)\varphi^{\prime 2}(Y)]\rightarrow\int_{-\frac{1}{2}}^{\frac{1}{2}}\theta^2(t)dt\mathbb{E}[\varphi^{\prime 2}(Y)]$ what proves the second one.\hfill$\diamond$\\

If $\mathbb{P}_Y$ is absolutely continuous and satisfies the Hamza condition, Fukushima {\it et al.} (1991, theorem 3.1.6 p.105), e.g. as soon as $\mathbb{P}_Y$ has a continuous density, the form $\mathcal{E}[\varphi]=\frac{1}{24}\mathbb{E}_Y[\varphi^{\prime 2}]$ is Dirichlet and admits the square field operator $\Gamma[\varphi]=\frac{1}{12}\varphi^{\prime 2}$. The graduation yields therefore an error structure $(\mathbb{R}, \mathcal{B}(\mathbb{R}), \mathbb{P}_Y,\mathbb{D},\Gamma)$, Bouleau (2003, Chapter III) whose operator $\Gamma$ does not depend on $\mathbb{P}_Y$ provided that $Y$ has a regular density. This translates in terms of errors the arbitrary functions principle.\\

\noindent{\bf Remark.} In his study of the exception to the law of errors of Gauss, Henri Poincar\'e (1912, p.219) considers the case where the nearest graduation is chosen except in a central zone across the middle between two marks of the graduation, where either the left one or the right one is chosen randomly with equal probability. It is easily seen that this case may be handled in the same way and gives similar results. The Dirichlet form is increased with respect to the above case.\\

\noindent{\large\textsf{I.3.}} {\bf Historical comment.}

In his intuitive version, the idea underlying the arbitrary functions method is ancient. The historian J. von Plato (1983) dates it back to a book of J. von Kries (1886). We find indeed in this philosophical treatise the idea that if a roulette had equal and infinitely small black and white cases, then there would be an equal probability to fall on a case or on the neighbour one, hence by addition an equal probability to fall either on black or on white. But no precise proof was given. The idea remains at the common sense level. 

A mathematical argument for the fairness of the roulette and for the equi-distribu\-tion of other mechanical systems (little planets on the Zodiac) was proposed by H. Poincar\'e  in his course on probability published in 1912 (1912, Chap. VIII \S 92 and especially \S 93). In present language, Poincar\'e shows the weak convergence of $tX+Y \mbox{mod }  2\pi$ when $t\uparrow \infty$ to the uniform law on $(0,2\pi)$ when the pair $(X,Y)$ has a density. He uses the characteristic functions. His proof supposes the density  be $\mathcal{C}^1$ with bounded derivative in order to perform an integration by parts, but the proof would extend to the general absolutely case if we were using instead the Riemann-Lebesgue lemma. 

The question is then developed  without major changes by several authors, E. Borel (1924) (case of continuous density), M. Fr\'echet (1921) (case of Riemann-integrable density), B. Hostinski (1926) (1931) (bidimensional case) and is tackled anew by E. Hopf (1934), (1936) and (1937) with the more general point of view of asymptotic behaviour of dissipative dynamical systems. Hopf has shown that these phenomena may be mathematically understood in the framework of ergodic theory and are related to mixing. Today the connection is close to Rajchman (or mixing) measures cf. Lyons (1995), Katok and Thouvenot (2005), interesting objects related to deep properties of descriptive set theory.\\

\newpage
\noindent{\Large\textsf{II. Theoretical tools}}\\

\noindent{\large\textsf{II.1.}}{\bf Approximation, Dirichlet forms and bias operators.}

Our study uses the theoretical framework concerning the bias operators and the Dirichlet form generated by an approximation proposed in Bouleau (2006). We recall here the definitions and main results for the convenience of the reader. Here, considered Dirichlet forms are always symmetric.

Let $Y$ be a random variable  defined on $(\Omega, \mathcal{A}, \mathbb{P})$ with values in a measurable space $(E,\mathcal{F})$ and let $Y_n$ be approximations also defined on $(\Omega, \mathcal{A}, \mathbb{P})$ with values in $(E,\mathcal{F})$.We consider an algebra $\mathcal{D}$ of bounded functions from $E$ into $\mathbb{R}$ or $\mathbb{C}$ containing the constants and dense in $L^2(E,\mathcal{F},\mathbb{P}_Y)$ and a sequence $\alpha_n$ of positive numbers. With $\mathcal{D}$ and $(\alpha_n)$ we consider the four following assumptions defining the four bias operators 
$$
(\mbox{H}1)\qquad\left\{\begin{array}{l}
\forall \varphi\in{\cal D}, \mbox{ there exists } \overline{A}[\varphi]\in L^2(E,{\cal F},\mathbb{P}_Y)\quad s.t. \quad\forall \chi\in{\cal D}\\
\lim_{n\rightarrow\infty} \alpha_n\mathbb{E}[(\varphi(Y_n)-\varphi(Y))\chi(Y)]=\mathbb{E}_Y[\overline{A}[\varphi]\chi].
\end{array}\right.
$$
$$
(\mbox{H}2)\qquad\left\{\begin{array}{l}
\forall \varphi\in{\cal D}, \mbox{ there exists } \underline{A}[\varphi]\in L^2(E,{\cal F},\mathbb{P}_Y)\quad s.t. \quad\forall \chi\in{\cal D}\\
\lim_{n\rightarrow\infty} \alpha_n\mathbb{E}[(\varphi(Y)-\varphi(Y_n))\chi(Y_n)]=\mathbb{E}_Y[\underline{A}[\varphi]\chi].
\end{array}\right.
$$
$$
(\mbox{H}3)\quad\left\{\begin{array}{l}
\forall \varphi\in{\cal D}, \mbox{ there exists } \widetilde{A}[\varphi]\in L^2(E,{\cal F},\mathbb{P}_Y)\quad s.t. \quad\forall \chi\in{\cal D}\\
\lim_{n\rightarrow\infty} \alpha_n\mathbb{E}[(\varphi(Y_n)-\varphi(Y))(\chi(Y_n)-\chi(Y))]=-2\mathbb{E}_Y[\widetilde{A}[\varphi]\chi].
\end{array}\right.
$$
$$
(\mbox{H}4)\quad\left\{\begin{array}{l}
\forall \varphi\in{\cal D}, \mbox{ there exists } \Abar[\varphi]\in L^2(E,{\cal F},\mathbb{P}_Y)\quad s.t. \quad\forall \chi\in{\cal D}\\
\lim_{n\rightarrow\infty} \alpha_n\mathbb{E}[(\varphi(Y_n)-\varphi(Y))(\chi(Y_n)+\chi(Y))]=2\mathbb{E}_Y[\Abar[\varphi]\chi].
\end{array}\right.
$$
We first note that as soon as two of hypotheses (H1) (H2) (H3) (H4) are fulfilled (with
 the same algebra ${\cal D}$ and the same sequence $\alpha_n$), the other two follow thanks to the relations
$$\widetilde{A}=\frac{\overline{A}+\underline{A}}{2}\quad\quad\Abar=\frac{\overline{A}-\underline{A}}{2}.$$
When defined, the operator $\overline{A}$ which considers the asymptotic error from the point of view of the limit model, will be called {\it the 
theoretical bias operator}.

The operator $\underline{A}$ which considers the asymptotic error from the point of view of the approximating model will be called {\it the 
practical bias operator}.

Because of the property
$$<\widetilde{A}[\varphi],\chi>_{L^2(\mathbb{P}_Y)}=<\varphi,\widetilde{A}[\chi]>_{L^2(\mathbb{P}_Y)}$$
the operator $\widetilde{A}$ will be called {\it the symmetric bias operator}.

The operator $\Abar$ which is often (see theorem 3 below) a first order operator will be called {\it the singular bias operator}.\\

\noindent{\bf Theorem 2.} {\it Under the hypothesis {\rm (H3)},

a) the limit
\begin{equation}\widetilde{\cal E}[\varphi,\chi]=\lim_n  \frac{\alpha_n}{2}\mathbb{E}[(\varphi(Y_n)-\varphi(Y))(\chi(Y_n)-\chi(Y)]\qquad \varphi, \chi\in{\cal D}\end{equation}
defines a closable positive bilinear form whose smallest closed extension is denoted $({\cal E},\mathbb{D})$.

b) $({\cal E},\mathbb{D})$ is a Dirichlet form

c) $({\cal E},\mathbb{D})$ admits a square field operator $\Gamma$ satisfying $\forall \varphi,\chi\in{\cal D}$
\begin{equation}
\Gamma[\varphi]=\widetilde{A}[\varphi^2]-2\varphi\widetilde{A}[\varphi\end{equation}
\begin{equation}\mathbb{E}_Y[\Gamma[\varphi]\chi]=\lim_n\alpha_n\mathbb{E}[(\varphi(Y_n)-\varphi(Y))^2(\chi(Y_n)+\chi(Y))/2]\end{equation}
\indent d) $({\cal E},\mathbb{D})$ is local if and only if $\forall \varphi\in{\cal D}$
\begin{equation}\lim_n \alpha_n\mathbb{E}[(\varphi(Y_n)-\varphi(Y))^4]=0\end{equation}
this condition is equivalent to
$\quad\exists\lambda>2\quad\lim_n\alpha_n\mathbb{E}[|\varphi(Y_n)-\varphi(Y)|^\lambda]=0.$

e) If the form $({\cal E},\mathbb{D})$  is local, then the} principle of asymptotic error calculus
{\it is valid on 
$\widetilde{\cal D}=\{F(f_1,\ldots,f_p)\;:\;f_i\in{\cal D},\;\;F\in{\cal C}^1(\mathbb{R}^p,\mathbb{R})\}$
i.e.\\

$
\lim_n\alpha_n\mathbb{E}[(F(f_1(Y_n),\ldots,f_p(Y_n))-F(f_1(Y),\ldots,f_p(Y))^2]$\hfill

\hfill$=\mathbb{E}_Y[\sum_{i,j=1}^p F^\prime_i(f_1,\ldots,f_p)F^\prime_j(f_1,\ldots,f_p)\Gamma[f_i,f_j]].$}\\

An operator $B$ from ${\cal D}$ into $L^2(\mathbb{P}_Y)$ will be said to be a {\it first order operator} if it satisfies
$$B[\varphi\chi]=B[\varphi]\chi+\varphi B[\chi]\qquad\forall\varphi,\chi\in{\cal D}$$
\vspace{-.1cm}

\noindent{\bf Theorem 3.} {\it Under} (H1) {\it to} (H4). {\it If there is a real number $p\geq 1$ s.t.
$$\lim_n\alpha_n\mathbb{E}[(\varphi(Y_n)-\varphi(Y))^2|\psi(Y_n)-\psi(Y)|^p]=0\quad\forall\varphi,\psi\in{\cal D}$$
then $\Abar$ is first order.}

In particular, if the Dirichlet form is local, by the d) of theorem 2, the operator $\Abar$ is first order.\\

\noindent{\bf Example 1.} In order to deepen the discussion begun above in the introduction about the conditional law of $Y_n$ given $Y=y$, we give below a simple example where both the conditional law of $Y_n$ given $Y=y$ and the conditional law of $Y$ given $Y_n=y$ are Dirac measures and where nevertheless the approximation gives rise to a non-zero Dirichlet form.

Let us consider the unit interval and the dyadic representation of real numbers. If $Y$ is uniformly distributed we may write 
$Y=\sum_{k=0}^\infty \frac{a_k}{2^{k+1}}$
where the $a_k$ are independent identically distributed with law $\frac{1}{2}\delta_0+\frac{1}{2}\delta_1$.

Let us approximate $Y$ by $Y_n=\sum_{k=0}^{n-1} \frac{a_k}{2^{k+1}}+\frac{1}{2}\sum_{k=n}^\infty \frac{a_k}{2^{k+1}}$. We see that $Y$ and $Y_n$ are deterministically linked :
$$Y_n=Y-\frac{1}{2}\frac{\{2^n Y\}}{2^n}\quad\quad Y=Y_n+\frac{1}{2}\frac{\{2^n Y\}}{2^n}.$$ Now, it is easily seen that on the algebra $\mathcal{D}=\mathcal{L}\{e^{2i\pi kx},k\in \mathbb{Z}\}$ we have
$$3.4^n\mathbb{E}[(\varphi(Y_n)-\varphi(Y))(\psi(Y_n)-\psi(Y))]\rightarrow\mathbb{E}[\varphi^\prime\overline{\psi^\prime}],$$ what gives in the real domain the Dirichlet form $\mathcal{E}[\varphi]=\frac{1}{2}\mathbb{E}[\varphi^{\prime 2}]$.\hfill$\diamond$\\

\noindent{\large\textsf{II.2.}} {\bf  Girsanov theorem for Dirichlet forms.}

An error structure is a probability space $(\Omega,\mathcal{A}, \mathbb{P})$ equipped with a local Dirichlet form with domain $\mathbb{D}$ dense in $L^2(\Omega,\mathcal{A}, \mathbb{P})$ admitting a square field operator $\Gamma$, see Bouleau (2003). We denote $\mathcal{D}A$ the domain of the associated generator.\\

\noindent{\bf Theorem 4.} {\it Let $(\Omega,\mathcal{A}, \mathbb{P}, \mathbb{D},\Gamma)$ be an error structure. Let be $f\in\mathbb{D}\cap L^\infty$ such that $f>0$, $\mathbb{E}f=1$, We put $\mathbb{P}_1=f.\mathbb{P}$.

a) The bilinear form $\mathcal{E}_1$ defined on $\mathcal{D}A\cap L^\infty$ by
\begin{equation}\mathcal{E}_1[u,v]=-\mathbb{E}\left[fvA[u]+\frac{1}{2}v\Gamma[u,f]\right]\end{equation}
is closable in $L^2(\mathbb{P}_1)$ and satisfies for $u,v\in\mathcal{D}A\cap L^\infty$
\begin{equation}\mathcal{E}_1[u,v]=-<A_1u,v>=-<u,A_1v>=\frac{1}{2}\mathbb{E}[f\Gamma[u,v]]\end{equation}
where $A_1[u]=A[u]+\frac{1}{2f}\Gamma[u,f]$.

b) Let $(\mathbb{D}_1,\mathcal{E}_1)$ be the smallest closed extension of $(\mathcal{D}A\cap L^\infty,\mathcal{E}_1)$. Then $\mathbb{D}\subset\mathbb{D}_1$, $\mathcal{E}_1$ is local and admits a square field operator $\Gamma_1$, and 
$$\Gamma_1=\Gamma\quad\mbox{\rm on}\quad \mathbb{D}$$
in addition $\mathcal{D}A\subset\mathcal{D}A_1$ and 
$A_1[u]=A[u]+\frac{1}{2f}\Gamma[u,f]$  for all $u\in\mathcal{D}A$.}\\

\noindent{\bf Proof.} 1) First, using that the resolvent operators are bounded operators sending $L^\infty$ into $\mathcal{D}A\cap L^\infty$, we see that $\mathcal{D}A\cap L^\infty$ is dense in $\mathbb{D}$ (equipped with the usual norm $(\|.\|_{L^2}^2+\mathcal{E}[.])^{1/2}$), hence also dense in $L^2(\mathbb{P}_1)$.

2) Using that $\mathbb{D}\cap L^\infty$ is an algebra, for $u,v\in\mathcal{D}A\cap L^\infty$ we have
$$\mathcal{E}_1[u,v]=-\mathbb{E}[fvA[u]+\frac{1}{2}v\Gamma[u,f]]=\frac{1}{2}\mathbb{E} [\Gamma[fv,u]-v\Gamma[u,f]]=\frac{1}{2}\mathbb{E}[f\Gamma[u,v]].$$
So, defining $A_1$ as in the statement, we have $\forall u,v\in\mathcal{D}A\cap L^\infty$
$$\mathcal{E}_1[u,v]=-\mathbb{E}_1[vA_1u]=-\mathbb{E}_1[uA_1v].$$
The operator $A_1$ is therefore symmetric on $\mathcal{D}A\cap L^\infty$ under $\mathbb{P}_1$. Hence the form $\mathcal{E}_1$ defined on $\mathcal{D}A\cap L^\infty$ is closable, Fukushima {\it et al.} (1994, condition 1.1.3 p 4).

3) Let $(\mathbb{D}_1,\mathcal{E}_1)$ be the smallest closed extension of $(\mathcal{D}A\cap L^\infty,\mathcal{E}_1)$. Let be $u\in\mathbb{D}$ and $u_n\in\mathcal{D}A\cap L^\infty$, with $u_n\rightarrow u$ in $\mathbb{D}$. Using $\mathcal{E}_1[u_n-u_m]\leq \|f\|_\infty\mathcal{E}[u_n-u_m]$ and the closedness of $\mathcal{E}_1$ we get $u_n\rightarrow u$ in $\mathbb{D}_1$, hence $\mathbb{D}\subset\mathbb{D}_1$. Now by usual inequalities we see that $\Gamma[u_n]$ is a Cauchy sequence in $L^1(\mathbb{P}_1)$ and that the limit $\Gamma_1[u]$ does not depend on the particular sequence $(u_n)$ satisfying the above condition. Then following Bouleau (2003, Chap. III \S2.5 p.38) the functional calculus extends to $\mathbb{D}_1$, the axioms of error structures are fulfilled for $(\Omega,\mathcal{A}, \mathbb{P}_1, \mathbb{D}_1,\Gamma_1)$ and this gives with usual arguments the b) of the statement.\hfill$\diamond$\\

\noindent{\large\textsf{II.3.}} {\bf Rajchman measures.}

\noindent{\bf Definition 1.} {\it A measure $\mu$ on the torus $\mathbb{T}^1$ is said to be Rajchman if 
$$\hat{\mu}=\int_{\mathbb{T}^1}e^{2i\pi nx}\,d\mu(x)\rightarrow 0\quad\quad\mbox{when }|n|\uparrow\infty.$$}

The set of Rajchman measures $\mathcal{R}$ is a band : if $\mu\in\mathcal{R}$ and if $\nu\ll|\mu|$ then $\nu\in\mathcal{R}$, cf.  Rajchman (1928) (1929), Lyons (1995).

\noindent{\bf Lemma.} {\it Let $X$ be a real random variable and let $\Psi_X(u)=\mathbb{E}e^{iuX}$ be its characteristic function. Then 
$$\lim_{|u|\rightarrow\infty}\Psi_X(u)=0\quad\Longleftrightarrow\quad\mathbb{P}_{\{X\}}\in\mathcal{R}.$$}

\noindent{\bf Proof.} a) If $\lim_{|u|\rightarrow\infty}\Psi_X(u)=0$ then  $\Psi_X(2\pi n)=
(\mathbb{P}_{\{X\}})\hat{}\,(n)\rightarrow 0.$

b) Let $\rho$ be a probability measure on $\mathbb{T}^1$ s.t. $\rho\in\mathcal{R}$. From
$$e^{2i\pi ux}=e^{2i\pi[u]x}\sum_{p=0}^\infty \frac{((u-[u])2i\pi x)^p}{p!}$$ we have
$$\int e^{2i\pi ux}\rho(dx)=\sum_{p=0}^\infty \frac{((u-[u])2i\pi )^p}{p!}\int x^p e^{2i\pi[u]x}\rho(dx).$$ Since $x^p\rho(dx)\in\mathcal{R}$ 
$$\int e^{2i\pi ux}\rho(dx)=\sum_{p=0}^\infty \frac{((u-[u])2i\pi )^p}{p!}a_p([u])$$with
$|a_p|\leq 1$ and $\lim_{|n|\rightarrow\infty}a_p(n)=0$, so 
$$\lim_{|u|\rightarrow\infty}\int e^{2i\pi ux}\rho(dx)=0.$$
Now if $\mathbb{P}_{\{X\}}\in\mathcal{R}$, since $1_{\{x\in[p,p+1[\}}.\mathbb{P}_{\{X\}}\ll\mathbb{P}_{\{X\}}$ we have
$$\lim_{|u|\rightarrow\infty}\mathbb{E}[e^{2i\pi uX}]=\lim_{|u|\rightarrow\infty}\sum_p \mathbb{E}[e^{2i\pi uX}1_{\{X\in[p,p+1[\}}]$$ which goes to zero by dominated convergence.\hfill$\diamond$\\

A probability measure on $\mathbb{R}$ satisfying the conditions of the lemma  will be called Rajchman. \\

\noindent{\bf Examples.} Thanks to the Riemann-Lebesgue lemma, absolutely continuous measures are in $\mathcal{R}$. 
It follows from the lemma that if a measure $\nu$ satisfies $\nu\star\cdots\star\nu\in\mathcal{R}$ then $\nu\in\mathcal{R}$.

Let $\beta\in]0,\frac{1}{2}[$, let $K_1$, be the closed set obtained by taking away from the unit interval the open interval of length $1-2\beta$ centered on $1/2$, $K_k$ the closed set obtained by iterating these suppressions homothetically on each segment of $K_{k-1}$ and let $\mu_k$ be the probability measure $\frac{1}{|K_k|}1_{K_k}.dx$ The continuous measure $\mu$ weak limit of the $\mu_k$ carried by the perfect set $\cap_k K_k$ is in $\mathcal{R}$ iff $1/\beta$ is not a Pisot number (a Pisot number is a root of a polynomial with entire coefficients and with coefficient of highest degree term equal to 1, irreducible over $\mathbb{Q}$, such that the other roots have a modulus $<1$) cf. Kahane and Salem (1953).\\

\noindent{\bf Proposition 1.} {\it Let $X,Y,Z$ be random variables with values in $\mathbb{R}$, $\mathbb{R}$, and $\mathbb{R}^m$ resp. Then 
\begin{equation}(\{nX+Y\},X,Y,Z)\quad\stackrel{d}{\Longrightarrow}\quad(U,X,Y,Z)\end{equation}
where $U$ is uniform on the unit interval independent of $(X,Y,Z)$, if and only if $\mathbb{P}_X$ is Rajchman.}\\

\noindent{\bf Proof.} If $\mu$ is a probability measure on $\mathbb{T}^1\times\mathbb{R}^m$, let us put 
$$\hat{\mu}(k,\zeta)=\int e^{2i\pi kx+<\zeta,y>}\mu(dx,dy),$$ then $\mu_n\stackrel{d}{\Longrightarrow}\mu$ iff $\hat{\mu}_n(k,\zeta)\rightarrow\hat{\mu}(k,\zeta)$ $\forall k\in\mathbb{Z}$, $\forall\zeta\in\mathbb{R}^m$.

a) If $\mathbb{P}_X\in\mathcal{R}$ 
$$\hat{\mathbb{P}}_{(\{nX+Y\},X,Y,Z)}(k,\zeta_1,\zeta_2,\zeta_3)=\mathbb{E}[\exp\{2i\pi k(nX+Y)+i\zeta_1X+i\zeta_2Y+i<\zeta_3,Z>\}]$$
$$=\int e^{2i\pi knx}f(x)\mathbb{P}_{\{X\}}(dx)$$
with $f(x)=\mathbb{E}[\exp\{2i\pi kY+i\zeta_1X+i\zeta_2Y+i<\zeta_3,Z>\}|\{X\}=x]$. The fact that $f.\mathbb{P}_{\{X\}}\in\mathcal{R}$ gives the result.

b) Conversely, taking $(k,\zeta_1,\zeta_2,\zeta_3)= (1,0,-2\pi, 0)$ gives  $\hat{\mathbb{P}}_{\{X\}}(n)\rightarrow 0$ i.e. $\mathbb{P}_X\in\mathcal{R}$.\hfill$\diamond$\\

The preceding definitions and properties extend to $\mathbb{T}^d$ : a measure $\mu$ on $\mathbb{T}^d$ is said to be in $\mathcal{R}$ if $\hat{\mu}(k)\rightarrow0$ as $k\rightarrow\infty$ in $\mathbb{Z}^d$. The set of measures in $\mathcal{R}$ is a band. If $X $ is $\mathbb{R}^d$-valued, $\lim_{|u|\rightarrow\infty}\mathbb{E}e^{i<u,X>}=0$ is equivalent to $\mathbb{P}_{\{X\}}\in\mathcal{R}$ where $\{x\}=(\{x_1\},\ldots,\{x_d\})$.\\

\noindent{\Large\textsf{III. Finite dimensional cases.}}\\

In this part we carry on deeper with the basic example (section {\textsf{I.2}}) in the finite dimensional case. \\

\noindent{\large\textsf{III.1.}} We suppose $Y$ is $\mathbb{R}^d$-valued, measured with an equidistant graduation corresponding to an orthonormal rectilinear coordinate system, and estimated to the nearest graduation component by component. Thus we put
$$Y_n=Y+\frac{1}{n}\theta(nY)$$ with $\theta(y)=(\frac{1}{2}-\{y_1\},\cdots,\frac{1}{2}-\{y_d\})$.\\

\noindent{\bf Theorem 5.} {\it a) If $\mathbb{P}_Y$ is Rajchman and if $X$ is $\mathbb{R}^m$-valued

\begin{equation}(X,n(Y_n-Y))\quad\stackrel{d}{\Longrightarrow}\quad(X,(V_1,\ldots,V_d))\end{equation}
where the $V_i$'s are independent identically distributed uniformly distributed on $(-\frac{1}{2},\frac{1}{2})$ and independent of $X$.

\noindent For all $\varphi\in\mathcal{C}^1\cap lip(\mathbb{R}^d)$
\begin{equation}
(X,n(\varphi(Y_n)-\varphi(Y)))\quad\stackrel{d}{\Longrightarrow}\quad(X,\sum_{i=1}^dV_i\varphi^\prime_i(Y))
\end{equation}
\begin{equation}
n^2\mathbb{E}[(\varphi(Y_n)-\varphi(Y))^2Ü|Y\!=\!y]\rightarrow\frac{1}{12}\sum_{i=1}^d\varphi^{\prime 2}_i(y)\qquad\mbox{ in }L^1(\mathbb{P}_Y)
\end{equation}
in particular

\begin{equation}
n^2\mathbb{E}[(\varphi(Y_n)-\varphi(Y))^2]\rightarrow\mathbb{E}_Y[\frac{1}{12}\sum_{i=1}^d\varphi^{\prime 2}_i(y)].
\end{equation}
\indent b) If $\varphi$ is of class $\mathcal{C}^2$, the conditional expectation $n^2\mathbb{E}[\varphi(Y_n)-\varphi(Y)|Y=y]$ possesses a version $n^2(\varphi(y+\frac{1}{n}\theta(ny))-\varphi(y))$ independent of the probability measure $\mathbb{P}$ which converges in the sense of distributions to the function $\frac{1}{24}\bigtriangleup\varphi$.

c) If $\mathbb{P}_Y\ll dy $ on $\mathbb{R}^d$, $\forall\psi\in L^1([0,1])$
\begin{equation}(X,\psi(n(Y_n-Y)))\quad\stackrel{d}{\Longrightarrow}(X,\psi(V)).\end{equation}

d) We consider the bias operators on the algebra $\mathcal{C}^2_b$ of bounded functions with bounded derivatives up to order 2 with the sequence $\alpha_n=n^2$. If $\mathbb{P}_Y\in\mathcal{R}$ and if one of the following condition is fulfilled

i) $\forall i=1,\ldots,d$ the partial derivative $\partial_i\mathbb{P}_Y$ in the sense of distributions is a measure $\ll\mathbb{P}_Y$ of the form $\rho_i\mathbb{P}_Y$ with $\rho_i\in L^2(\mathbb{P}_Y)$

ii) $\mathbb{P}_Y=h1_{G}\frac{dy}{|G|}$  with $G$ open set, $h\in H^1\cap L^\infty(G)$, $h>0$

\noindent then hypotheses {\rm (H1)} to {\rm (H4)} are satisfied and 
$$\begin{array}{rl}
\overline{A}[\varphi]&=\frac{1}{24}\bigtriangleup\varphi\\
\widetilde{A}[\varphi]&=\frac{1}{24}\bigtriangleup\varphi+\frac{1}{24}\sum\varphi^\prime_i\rho_i\qquad\mbox {case i)}\\
\widetilde{A}[\varphi]&=\frac{1}{24}\bigtriangleup\varphi+\frac{1}{24}\frac{1}{h}\sum h^\prime_i\varphi^\prime_i\qquad\mbox {case ii)}\\
\Gamma[\varphi]&=\frac{1}{12}\sum \varphi^{\prime 2}_i.
\end{array}
$$}

\noindent{\bf Proof.}  The argument for relation (10) is similar to one dimensional case. The relation (11) comes from the Taylor expansion 
$\varphi(Y_n)-\varphi(Y)=$

$=\sum_{i=1}^d (Y_{n,i}-Y_i)\int_0^1\varphi^\prime_i(Y_{n,1},\ldots,Y_{n,i-1},Y_i+t(Y_{n,i}-Y_i),Y_{i+1},\ldots,Y_d)\,dt$

\noindent and the convergence
$$(X,\sum_i\theta(nY_i)\varphi^\prime_i(Y))\quad\stackrel{d}{\Longrightarrow}\quad(X,\sum_i\varphi^\prime_i(Y)V_i)$$
thanks to (10) and the following approximation in $L^1$ 
$$\mathbb{E}\left|\sum_i\theta(nY_i)\varphi^\prime_i(Y)-\sum_i\theta(nY_i)\int_0^1\varphi^\prime_i(\ldots,Y_i+t(Y_{n,i}-Y_i),\ldots)dt\right|\rightarrow0.$$
To prove the formulas (12) and (13) let us remark that

\begin{flushleft}
$n^2\mathbb{E}[(\varphi(Y_n)-\varphi(Y)^2|Y=y]=$
\end{flushleft}
$$=\mathbb{E}\left[\left|\sum_i\theta(nY_i)\int_0^1\varphi^\prime_i(\ldots,Y_i+t(Y_{n,i}-Y_i),\ldots)dt\right|^2|Y=y\right]$$
$$=\left|\sum_i\theta(ny_i)\int_0^1\varphi^\prime_i(y_1+\frac{1}{n}\theta(ny_1),\ldots,y_i+t\frac{1}{n}\theta(ny_i),\ldots)dt\right|^2\quad\mathbb{P}_Y-a.s.$$
each term $(\theta(ny_i)\int_0^1\varphi^\prime_i(\ldots)dt)^2$ converges to $\int\theta^2\varphi^{\prime 2}_i(y)=\frac{1}{12}\varphi^{\prime 2}_i$ in $L^1$ and each term $\theta(ny_i)\theta(ny_j)\int_0^1\ldots\int_0^1\ldots$ goes to zero in $L^1$ what proves the a) of the statement.

The point b) is obtained following the same lines with a Taylor expansion up to second order and an integration by part thanks to the fact that $\varphi$ is now supposed to be $\mathcal{C}^2$.

In order to prove c) let us suppose first that $\mathbb{P}_Y=1_{[0,1]^d}.dy$. Considering a sequence of functions $\psi_k\in\mathcal{C}_b$ tending to $\psi$ in $L^1$ we have the bound 
$$
\begin{array}{l}
|\mathbb{E}[e^{i<u,X>}e^{iv\psi(\theta(nY))}]-\mathbb{E}[e^{i<u,X>}e^{iv\psi_k(\theta(nY))}]|\\
\leq |v|\int|\psi(\theta(ny))-\psi_k(\theta(ny))|dy\\
=|v|\sum_{p_1=0}^{n-1}\cdots\int_{p_1}^{p_1+1}\cdots|\psi(\theta(ny_1)\ldots)-\psi_k(\theta(ny_1)\ldots)|dy_1\ldots dy_d\\
=|v|\sum\cdots\sum\int\cdots\int|\psi(\theta(x_1),\ldots)-
\psi_k(\theta(x_1),\ldots)|\frac{dx_1}{n}\cdots\frac{dx_d}{n}\\
=|v|\|\psi-\psi_k\|_{L^1}.
\end{array}
$$
This yields (14) in this case. Now if $\mathbb{P}_Y\ll dy$ then $\mathbb{P}_{\{Y\}}\ll dy $ on $[0,1]^d$ and the weak convergence under $dy$ on $[0,1]^d$ implies the weak convergence under $\mathbb{P}_{\{Y\}}$ what yields the result.

In d) the point i) is proved by the approach already used in Bouleau (2006) consisting of proving  that hypothesis (H3) is fulfilled by displaying the operator $\widetilde{A}$ thanks to an integration by parts. The point ii) is an application of Girsanov theorem (theorem 4).\hfill$\diamond$\\

\noindent{\bf Remarks.} 1) About the relations (11) (12) (13), let us note that with respect to the form $$\mathcal{E}[\varphi]=\frac{1}{24}\mathbb{E}_Y \sum_i\varphi^{\prime 2}_i$$ when it is closable, the random variable $\sum_iV_i\varphi^\prime_i$ appears to be {\it a gradient} :
if we put $\varphi^\#=\sum_iV_i\varphi^\prime_i$ then a we have
$$\mathbb{E}[\varphi^{\# 2}]=\frac{1}{12}\sum_i\varphi^{\prime 2}_i=\Gamma[\varphi]$$
the square field operator associated to $\mathcal{E}$. We will find this phenomenon again on the Wiener space.

2) If $d=1$, the convergence in (12) holds in $L^p$ $1\leq p<\infty$ and in d) belonging to $\mathcal{R}$ is automatic under i) or ii).

3) Let us also remark that when $d=1$ assumptions i) and ii) may be replaced by the Hamza condition on $\mathbb{P}_Y$ which suffices to imply (H1) to (H4).\\

\noindent{\large\textsf{III.2.}} {\bf Approximation to the nearest graduation, by excess, or by default.}

Let us come back to the basic example. When the approximation is done to the nearest graduation, on the algebra $\mathcal{C}^2_b$ the four bias operators are zero with the sequence $\alpha_n=n$ (cf. theorem 5 where $\alpha_n=n^2$).

But we would obtain a quite different result with an approximation by default or by excess where the effect of the shift is dominating.

If the random variable $Y$ is approximated by default by $Y_n^{(d)}=\frac{[nY]}{n}$ then 
$$n(Y_n^{(d)}-Y)\stackrel{d}{\Longrightarrow} -U\quad\mbox{and}\quad\mathbb{E}[n(Y_n^{(d)}-Y)]\rightarrow -\frac{1}{2}$$
as soon as $Y$ is say bounded. With this approximation, if we do not erase the shift down proportional to $-\frac{1}{2n}$, and if we take $\alpha_n=n$ we obtain  first order bias operators without diffusion :
$\overline{A}[\varphi]=-\frac{1}{2}\varphi^\prime=-\underline{A}[\varphi]$ and $\widetilde{A}=0$. The same happens of course with the approximation by excess.\\

\noindent{\large\textsf{III.3.}} {\bf Extension to more general graduations.}

Let $Y$ be an $\mathbb{R}^d$-valued random variable approximated by $Y_n=Y+\xi_n(Y)$ with a sequence $\alpha_n\uparrow\infty$ on the algebra $\mathcal{D}=\mathcal{L}\{e^{<u,x>},\;u\in\mathbb{R}^d\}$, the function $\xi_n$ satisfying 
$$
(\ast)\left\{\begin{array}{l}
\alpha_n\mathbb{E}[|\xi_n|^3(Y)]\rightarrow 0\\
\\
\alpha_n\mathbb{E}[\varphi(Y)<u,\xi_n(Y)>^2]\rightarrow\mathbb{E}_Y[\varphi . u^\ast\underline{\underline{\gamma}}u]\qquad\forall\varphi\in\mathcal{D},\forall u\in\mathbb{R}^d\\
\mbox{ with }\gamma_{ij}\in L^\infty(\mathbb{P}_Y)\mbox{ and }\frac{\partial \gamma_{ij}}{\partial x_j}\mbox{ in distributions sense }\in L^2(\mathbb{P}_Y)\\
\\
\alpha_n\mathbb{E}[\varphi(Y)<u,\xi_n(Y)>]\rightarrow 0\quad \forall\varphi\in\mathcal{D}.
\end{array}
\right.
$$
Under these hypotheses we have\\

\noindent{\bf Theorem 6.} {\it a) {\rm (H1)} is satisfied and 
$$\overline{A}[\varphi]=\frac{1}{2}\sum_{ij}\gamma_{ij}\frac{\partial^2\varphi}{\partial x_i\partial x_j}.$$
b) If for $i=1,\ldots,d$, the partial derivative $\partial_i\mathbb{P}_Y$ in the sense of distributions is a bounded measure of the form $\rho_i\mathbb{P}_Y$ with $\rho_i\in L^2(\mathbb{P}_Y)$ then assumptions {\rm (H1)} to {\rm (H4)} are fulfilled and $\forall\varphi\in\mathcal{D}$
$$\widetilde{A}[\varphi]=\frac{1}{2}\sum_{ij}\gamma_{ij}\frac{\partial^2\varphi}{\partial x_i\partial x_j}+\sum_i(\sum_j(\frac{\partial \gamma_{ij}}{\partial x_j}+\gamma_{ij}\rho_j))\frac{\partial \varphi}{\partial x_i}$$
the square field operator is 
$$\Gamma[\varphi]=\sum_{ij}\frac{\partial \varphi}{\partial x_i}\frac{\partial \varphi}{\partial x_j}.$$}\\

\noindent{\bf Proof.} By the choice of the algebra $\mathcal{D}$ this theorem is simple. The argument consists of elementary Taylor expansions in order to prove the existence of the bias operators. Then theorem 2 applies.\hfill$\diamond$\\

\noindent{\Large\textsf{IV. Rajchman martingales.}}\\

Let $(\mathcal{F}_t)$ be a right continuous filtration on $(\Omega, \mathcal{A}, \mathbb{P})$ and $M$ be a continuous local $(\mathcal{F}_t,\mathbb{P})$-martingale nought at zero. $M$ will be said to be Rajchman if the measure $d\langle M,M\rangle_s$ restricted to compact intervals belongs to $\mathcal{R}$ almost surely.\\

\noindent{\large\textsf{IV.1.}} We will show that the method followed by Rootz\'en (1980) extends to Rajchman martingales and provides the following\\

\noindent{\bf Theorem 7.} {\it Let $M$ be a continuous local martingale which is Rajchman and s.t. $\langle M,M\rangle_\infty=\infty$.

Let $f$ be a bounded Riemann-integrable periodic function with unit period on $\mathbb{R}$ s.t. $\int_0^1f(s)ds=0$.
We put $T_n(t)=\inf\{s\,:\,\int_0^sf^2(nu)\,d\langle M,M\rangle_u\;>\;t\}$. Then for any random variable $X$
\begin{equation}(X,\int_0^.f(ns)\,dM_s)\quad\stackrel{d}{\Longrightarrow}\quad(X,W_{\|f\|^2\langle M,M\rangle_.}),\end{equation}
the weak convergence is understood on $\mathbb{R}\times\mathcal{C}([0,1])$ and  $W$ is an independent standard Brownian motion.}\\

Before proving the theorem, let us remark that it shows that the random measure $dM_s$ behaves in some sense like a Rajchman measure. Indeed if $\mathbb{P}_Y\in\mathcal{R}$ we have 
$$\int_{-\infty}^yg(nx)\mathbb{P}_Y(dx)\rightarrow \int_0^1g(x)dx\int_{-\infty}^y\mathbb{P}_Y(dx)$$
as soon as $g$ is periodic with unit period, Riemann-integrable and bounded. Now applying the theorem to the Brownian motion gives the similar relation
$$\int_0^t f(ns)\,dB_s\quad\stackrel{d}{\Longrightarrow}\quad(\int_0^1f^2(s)ds)^{1/2}\int_0^tdW_s.$$
{\bf Proof.} We consider the local martingale $N_t=\int_0^tf(ns)dM_s$.

a) In order to be sure that $\langle N,N\rangle_\infty=\infty$, we change $N_t$ into $\tilde{N_t}=\int_0^tf_n(s)dM_s$ with $f_n(s)=f(ns)$ for $s\in[0,1)$, $f_n(s)=0$ for $s\in[1,n]$ and $f_n(s)=1$ for $t>n$. We put $S_n(t)=\inf\{s\,:\langle\tilde{N},\tilde{N}\rangle_s\;>\;t\}$.

b) We want to show
\begin{equation}\qquad\mathbb{E}[\xi F(\tilde{N}_{S_n})]\rightarrow\mathbb{E}[\xi F(W)]\quad\forall\xi\in L^1(\mathbb{P})\quad\forall F\in \mathcal{C}_b([0,1]).\end{equation}
It is enough to consider the case $\xi>0$, $\mathbb{E}\xi=1$, and $\xi$ may be supposed to be $\mathcal{F}_T$-measurable for a deterministic time $T$ large enough. Let be $\tilde{\mathbb{P}}=\xi.\mathbb{P}$ and $D(t)=\mathbb{E}[\xi|\mathcal{F}_t]$. The process
$$\tilde{M_t}=M_t-\int_0^t D^{-1}(s)d\langle M,D^c\rangle_s$$ is a continuous local martingale under $\tilde{\mathbb{P}}$.  Therefore $\int_0^{S_n(t)}f_n(s)\,d\tilde{M}_s$ is a Brownian motion under $\tilde{\mathbb{P}}$ Revuz and Yor (1994,  p.313 theorem 1.4 and p 173). Writing
$$\int_0^{S_n(t)}f_n(s)dM_s=\int_0^{S_n(t)}f_n(s)d\tilde{M}_s+\int_0^{S_n(t)}\frac{f_n(s)}{D(s)}d\langle M,D^c\rangle_s$$ and noting that $d\langle M,D^c\rangle_s$ vanishes on $]T,\infty[$, in order to show $(16)$ it suffices to show
$$\sup_{0\leq t\leq T}\left|\int_0^t\frac{f_n(s)}{D(s)}\,d\langle M,D^c\rangle_s\right|\rightarrow 0\quad \mbox{ a.s. when } n\rightarrow\infty$$
hence to show
$$\sup_{0\leq t\leq 1}\left|\int_0^t\frac{f(ns)}{D(s)}\,d\langle M,D^c\rangle_s\right|\rightarrow 0\quad \mbox{ a.s. when } n\rightarrow\infty$$
and, because $M$ is Rajchman this comes from the following lemma :\\

\noindent{\bf Lemma.} {\it Let $f$ be as in the statement of the theorem, then $\forall\mu\in\mathcal{R}$
$$\sup_{0\leq t\leq 1}\left|\int_0^tf(ns)\mu(ds)\right|\rightarrow 0\quad \mbox{ as } n\rightarrow\infty.$$}

\noindent{\bf Proof.} We have 
$$\int_0^t f(ns)\mu(ds)\rightarrow\int_0^1f(s)ds\int_0^t\mu(ds)=0.$$
Since $f$ is bounded, the functions $\int_0^tf(ns)\mu(ds)$ are equi-continuous and the result follows from Ascoli theorem.\hfill$\diamond$\\

c) This proves the following stable convergence
$$(X,\int_0^{T_n(.)}f(ns)dM_s)\quad\stackrel{d}{\Longrightarrow}\quad(X,W_.)$$ and by the fact that the following limit 
$$\int_0^t f^2(ns)d\langle M,M\rangle_s\rightarrow\int_0^1f^2(s)ds\langle M,M\rangle_t$$ is a continuous process, this gives the announced result.\hfill$\diamond$\\

\noindent{\bf Remark.} If $\int_0^1f(s)ds\neq 0$, then keeping the other hypotheses unchanged, we obtain
$$(X,\int_0^.f(ns)dM_s)\quad\stackrel{d}{\Longrightarrow}\quad\left(X,(\int_0^1f(s)ds)M_.+(\int_0^1(f-\!\int_0^1f)^2)^{1/2}W_{\langle M,M\rangle_.}\right).$$
\vspace{0.2cm}

\noindent{\large\textsf{IV.2.}} {\bf Limit quadratic form for Rajchman martingales.}\\

We study the induced limit quadratic form when the martingale $M$ is approximated by the martingale $M^n_t=M_t+\int_0^t\frac{1}{n}f(ns)dM_s$. The notation is the same as in the preceding section and $f$ satisfies the same hypotheses as in theorem 7.\\

\noindent{\bf Theorem 8.} {\it Let $M$ be a Rajchman martingale s.t. $M_1\in L^2$ and $\eta$, $\zeta$ bounded adapted processes. Then
$$
\begin{array}{c}
\displaystyle n^2\mathbb{E}\left[(\exp\{i\int_0^1\eta_sdM^n_s\}-\exp\{i\int_0^1\eta_sdM_s\})(\exp\{i\int_0^1\zeta_sdM^n_s\}-\exp\{i\int_0^1\zeta_sdM_s\})\right]\\
\\
\displaystyle\rightarrow-\mathbb{E}\left[\exp\{i\int_0^1(\eta_s+\zeta_s)dM_s\}\int_0^1\eta_s\zeta_s\,d\langle M,M\rangle_s\right]\int_0^1f^2(s)ds.
\end{array}
$$}

{\bf Proof.} By the finite increments formula, the first term in the statement may be written
$$-\mathbb{E}[\exp\{i\int_0^1(\eta_s+\zeta_s)dM_s\}\int_0^1\eta_s f(ns)dM_s\int_0^1\zeta_sf(ns)dM_s]+o(1)$$ therefore, thanks to theorem 7, the theorem is a consequence of the following lemma :\\

\noindent{\bf Lemma.} {\it Suppose $\mathbb{E}M_1^2<\infty$ and $\eta$ adapted and bounded, then the random variables $\int_0^1\eta_sf(ns)dM_s$ are uniformly integrable.}\\

\noindent{\bf Proof.} It suffices to remark that their $L^2$-norm is equal to $\mathbb{E}\int_0^1\eta_s^2f^2(ns)\,d\langle M,M\rangle_s$ hence uniformly bounded.\hfill$\diamond$\\

\noindent{\Large\textsf{V. Sufficient closability conditions on the Wiener space.}}\\

The closability problem of the limit quadratic forms obtained in the preceding section, may be tackled with the tools available on the Wiener space.

Let us approximate the Brownian motion $(B_t)_{t\in[0,1]}$ by the process $B^n_t=B_t+\int_0^t\frac{1}{n}f(ns)\,dB_s$ where $f$ satisfies the same hypotheses as before. We  consider here only deterministic integrands.\\

\noindent{\bf Theorem 9.} {\it a) Let $\xi\in L^2([0,1])$, and let $X$ be a random variable defined on the Wiener space, i.e. a Wiener functional, then
\begin{equation}
\left(X,n(\exp\{i\int_0^1\xi dB^n\}-\exp\{i\int_0^1\xi dB\})\right)
\stackrel{d}{\Longrightarrow}\left(X,\|f\|_{L^2}(\exp\{i\int_0^1\xi dB\})^{\#}\right)
\end{equation}
here for any regular Wiener functional $Z$ we put $Z^\#(\omega,w)=\int_0^1D_sZ\,dW_s$, where $W$ is an independent Brownian motion.

b) 
\begin{equation}
n^2\mathbb{E}\left[(e^{i\xi.B^n}-e^{i\xi.B})^2\right]\rightarrow-\mathbb{E}[e^{2i\xi.B}]\int_0^1\xi^2ds\|f\|^2_{L^2}\end{equation}
on the algebra $\mathcal{L}\{e^{i\xi.B}\}$ the quadratic form $-\frac{1}{2}\mathbb{E}[e^{2i\xi.B}]\int_0^1\xi^2ds$ is closable, its closure is the Ornstein-Uhlenbeck form.}\\

\noindent{\bf Proof.} a) The first assertion comes easily from the similar result concerning Rajchman martingales using the fact that
$\int_0^1e^{i\alpha\int_0^1\frac{1}{n}f(ns)dB_s}d\alpha\rightarrow1$ in $L^p$ $p\in[1,\infty[$.

b) The obtained quadratic form is immediately recognized as the Ornstein-Uhlenbeck form which is closed. It follows that hypothesis (H3) is fulfilled and the symmetric bias operator is
$$\widetilde{A}[e^{i\int\xi dB}]=\left(-\frac{i}{2}\int\xi dB-\frac{1}{2}\int\xi^2 ds\right)e^{i\int \xi dB}.$$\hfill$\diamond$\\

If instead of the Wiener measure $m$, we consider the measure $m_1=h.m$ for an $h>0$, $h\in\mathbb{D}_{ou}\cap L^\infty$ where $\mathbb{D}_{ou}$ $(=D^{2,1})$ denotes the domain of the Ornstein-Uhlenbeck form, we know by Girsanov theorem (theorem 4) that the form 
$-\frac{1}{2}\mathbb{E}_1[e^{2i\xi.B}\int_0^1\xi^2ds]$ is closable, admits the same square field operator on $\mathbb{D}_{ou}$, and that its generator $A_1$ satisfies
$$A_1[\varphi]=\widetilde{A}[\varphi]+\frac{1}{2h}\Gamma_{ou}[\varphi,h]\quad\mbox{ for }\quad\varphi\in\mathcal{D}A_{ou}$$
Since the point a) of the theorem is still valid under $m_1$ because of the properties of stable convergence, the preceding theorem is still valid under $m_1$, the Dirichlet form being now
$$\mathcal{E}_1[\varphi]=\frac{1}{2}\mathbb{E}_1[\Gamma_{ou}[\varphi]]\quad\mbox{ for }\quad \varphi\in\mathbb{D}_{ou}.$$ 

\noindent{\bf Remark.} Let us come back to the general case of Rajchman martingales. If we suppose the Rajchman local martingale $M$ is in addition Gaussian, which is equivalent to suppose $\langle M,M\rangle$ deterministic, then  on the algebra $\mathcal{L}\{e^{i\int\xi dM}; \xi$  deterministic bounded $\}$ the limit quadratic form 
$$-\mathbb{E}[e^{i\int(\eta+\zeta)dM}\int_0^1\zeta_s\eta_sd\langle M,M\rangle_s]\|f\|^2_{L^2}$$
is closable, hence (H3) is satisfied.

Indeed, it suffices to exhibit the corresponding symmetric bias operator. But by the use of the calculus for Gaussian variables, it is easily seen that the operator defined by
$$\widetilde{A}[e^{i\int\xi dM}]=e^{i\int\xi dM}\left(-\frac{i}{2}\int\xi dM-\frac{1}{2}\int\xi^2\,d\langle M,M\rangle_s\right)\int f^2 ds$$ satisfies the required condition.$\hfill$$\diamond$\\

In the case $f(x)=\theta(x)=\frac{1}{2}-\{x\}$, the approximation used in the parts {\textsf{IV}} and {\textsf{V}} consists in approximating $B_t$ by $B_t+\int_0^t\frac{1}{n}\theta(ns)dB_s$. It is the most natural approximation suggested by the Rajchman property and the arbitrary functions principle. It yields also other approximation operators on the Wiener space, cf. Bouleau (CRAS 2006).

But it is different from the approximations usually encountered in the discretization of stochastic differential equations.

In order to draw a  link between the preceding study and works concerning the discretization of stochastic differential equations by the Euler scheme, especially those of Kurtz and Protter (1991), and Jacod and Protter (1998),  we may remark that the preceding results which yield 
\begin{equation}
\left(n\int_0^.(s-\frac{[ns]}{s})dB_s,n\int_0^.(B_s-B_{\frac{[ns]}{n}})ds,B_.\right)\stackrel{d}{\Longrightarrow}(\frac{1}{\sqrt{12}}W_.+\frac{1}{2}B_.,-\frac{1}{\sqrt{12}}W_.+\frac{1}{2}B_.,B_.)
\end{equation}
are generally hidden by a dominating phenomenon
$$
\left(\sqrt{n}\int_0^.(B_s-B_{\frac{[ns]}{n}})dB_s,B_.\right)\stackrel{d}{\Longrightarrow}(\frac{1}{\sqrt{2}}\widetilde{W_.},,B_.)
$$
due to the fact that when a variable of the second chaos (or in  further chaos) converges stably to a Gaussian limit, this one appears to be independent of the first chaos and therefore of $B$ itself.

The stable convergence (19) acts even on the first chaos. It concerns, for example, stochastic differential equations of the form
\begin{equation}\left\{
\begin{array}{l}
X^1_t=x^1_0+\int_0^tf^{11}(X^2_s)dB_s+\int_0^tf^{12}(X^1_s,X^2_s)ds\\
X^2_t=x^2_0+\int_0^tf^{22}(X^1_s,X^2_s)ds
\end{array}\right.
\end{equation}
where $X^1$ is with values in $\mathbb{R}^{k_1}$,  $X^2$ in $\mathbb{R}^{k_2}$,  $B$ in $\mathbb{R}^d$ and $f^{ij}$ are matrices with suitable dimensions. Such equations are encountered  to describe the movement of mechanical systems under the action of forces with a random noise, when the noisy forces depend only on the position of the system and the time. Typically 
$$\left\{
\begin{array}{l}
X_t=X_0+\int_0^tV_sds\\
V_t=V_0+\int_0^ta(X_s,V_s,s)ds+\int_0^tb(X_s,s)dB_s
\end{array}\right.
$$
which is a perturbation of the equation $\frac{d^2x}{dt^2}=a(x,\frac{dx}{dt},t)$. In such equations the stochastic integral may be understood as Ito as well as Stratonovitch. For the equation (20) the method of Kurtz and Protter (1991) without major changes yields the following result that we state in the case $k_1=k_2=d=1$ for simplicity.\\

\noindent{\bf Theorem 10.} {\it If functions $f^{ij}$ are $\mathcal{C}^1_b$, and if $X^n$ is the solution of {\rm(22)} by the Euler scheme, 
$$(n(X^n-X),X,B)\stackrel{d}{\Longrightarrow}(U,X,B)
$$
where the process $U$ is solution of the stochastic differential equation
$$
U"(t)=\sum_{k,j}\int_0^t\frac{\partial f^{ij}}{\partial x_k}(X_s)U^k_sdY^j_s-\sum_{k,j}\int_0^t\frac{\partial f^{ij}}{\partial x_k}(X_s)\sum_m f^{km}(X_s)dZ^{mj}_s
$$ where $Y_s=(B_s,s)^t $ and 
$$\begin{array}{l}
dZ^{12}_s=\frac{1}{\sqrt{12}}dW_s+\frac{1}{2}dB_s\\
dZ^{21}_s=-\frac{1}{\sqrt{12}}dW_s+\frac{1}{2}dB_s\\
dZ^{22}_s=\frac{ds}{2}
\end{array}
$$  
and as ever $W$ is an independent Brownian motion.}

\begin{list}{}
{\setlength{\itemsep}{0cm}\setlength{\leftmargin}{0.5cm}\setlength{\parsep}{0cm}\setlength{\listparindent}{-0.5cm}}
  \item\begin{center}
{\small REFERENCES}
\end{center}\vspace{0.4cm}
 {\sc Borel E.} {\it Calcul des probabilit\'es} Paris, 1924.

 {\sc Bouleau N.} {\it Error Calculus for Finance and Physics, the Language of Dirichlet Forms}, De Gruyter, 2003.

{\sc Bouleau N.}  ``When and how an error yields a Dirichlet form" {\it J. Funct. Anal.}, Vol 240, n$^0$2, (2006), 445-494.. 

{\sc Bouleau N.}  ``An extension to the Wiener space of the arbitrary functions principle" C. R. Acad. Sci. Paris,Ser I 343 (2006) 329-332.

{\sc Fr\'echet M.} ``Remarque sur les probabilit\'es continues"{\it Bull. Sci. Math.} $2^e$ s\'erie, 45, (1921), 87-88.

 {\sc Fukushima, M.; Oshima, Y.; Takeda, M.} {\it Dirichlet forms and symmetric Markov processes}, De Gruyter 1994.

 {\sc Hopf E.} ``On causality, statistics and probability" {\it J. of Math. and Physics} 18 (1934) 51-102.

{\sc Hopf E.} ``\"{U}ber die Bedeutung der willk\"{u}rlichen Funktionen f\"{u}r die Wahrscheinlichkeitstheorie" {\it Jahresbericht der Deutschen Math. Vereinigung} XLVI, I, 9/12, 179-194, (1936).

 {\sc Hopf E.} ``Ein Verteilungsproblem bei dissipativen dynamischen Systemen" {\it Math. Ann.} 114, (1937), 161-186.

 {\sc Hostinsk\'y B.} ``Sur la m\'ethode des fonctions arbitraires dans le calcul des probabilit\'es" {\it Acta Math.} 49, (1926), 95-113.

 {\sc Hostinsk\'y B.} {\it M\'ethodes g\'en\'erales de Calcul des Probabilit\'es}, Gauthier-Villars 1931.

 {\sc Jacod, J., Protter, Ph.} ``Asymptotic error distributions for the Euler method for stochastic 
differential equations'' {\it Ann. Probab.} 26, 267-307, (1998)

{\sc Katok A., Thouvenot J.-P.} ``Spectral methods and combinatorial constructions in ergodic theory" {\it Handbook in Dynamical systems} vol 1B, Elsevier, (2005), 649-743.

 {\sc von Kries, J.} {\it Die Prinzipien der Wahrscheinlichkeitsrechnung}, Freiburg 1886.

{\sc Kahane j.-P., Salem R.} {\it Ensembles parfaits et s\'eries trigonom\'etriques}, Hermann (1963).

 {\sc Kurtz, Th; Protter, Ph.} ``Wong-Zakai corrections, random evolutions and simulation schemes for SDEs" {\it Stochastic Analysis} 331-346,
Acad. Press, 1991.

 {\sc Lyons, R.}``Seventy years of Rajchman measures" {\it J. Fourier Anal. Appl.} Kahane special issue (1995), 363-377.

 {\sc von Plato, J.} ``The Method of Arbitrary Functions" {\it Brit. J. Phil. Sci.} 34, (1983), 37-42.

{\sc Poincar\'e, H.} {\it Calcul des Probabilit\'es} Gauthier-Villars, 1912.

{\sc Rajchman A.} ``Sur une classe de fonctions \`a variation born\'ee" {\it C. R. Acad. Sci. Paris}187, (1928), 1026-1028.

{\sc Rajchman A.}``Une classe de s\'eries g\'eom\'etriques qui convergent presque partout vers z\'ero" {\it Math. Ann.}101, (1929), 686-700.

{\sc Revuz D. Yor M.} {\it Continuous Martingales and Brownian motion}, Springer 1994.

 {\sc Rootz\'en, H.} ``Limit distribution for the error in approximation of stochastic integrals'' {\it Ann. Prob. } 8, 241-251, (1980).

\end{list}

\end{document}